\documentclass[a4paper,12pt]{amsart}
\usepackage{amssymb}
\usepackage{graphicx}
\usepackage{hyperref}
\usepackage{amsmath}
\usepackage{amsthm}
\usepackage{amssymb}
 
\hypersetup{colorlinks =true,citecolor=blue,linkcolor=blue,filecolor=blue,urlcolor=blue,citecolor=blue}
\usepackage{ifthen}
\nonstopmode \numberwithin{equation}{section}
\newtheorem{thm}{Theorem}%[section]
\newtheorem{cor}{Corollary}%[section]
\newtheorem{lem}{Lemma}%[section]
\newtheorem{conj}{Conjecture}

\theoremstyle{definition}

\newtheorem{prob}[equation]{Problem}

\newenvironment{rem}{%
	\bigskip
	\noindent \textsl{{\sl Remark. }}}{\bigskip}
\newenvironment{rems}{%
	\bigskip
	\noindent \textsl{{\sl Remarks. }}}{\bigskip}

\newcounter {own}
\def\theown {\thesection       .\arabic{own}}

\newenvironment{pf}[1][]{%
	\vskip 3mm
	\noindent
	\ifthenelse{\equal{#1}{}}%
	{{\slshape Proof. }}%
	{{\slshape #1.} }%
}%
{\qed\bigskip}

\newcounter{alphabet}

\newcommand{\IR}{{\mathbb R}}
\newcommand{\R}{{\mathbb R}}

\newcommand{\IC}{{\mathbb C}}
\newcommand{\C}{{\mathbb C}}
\newcommand{\IA}{{{\mathbb I}_1}}
\renewcommand{\Re}{{\,\operatorname{Re}\,}}
\renewcommand{\Im}{{\,\operatorname{Im}\,}}
\newcommand{\Ip}{{\mathbb{I}^{+}}}

\newcommand{\uhp}{{\mathbb H}}
\newcommand{\Hm}{{\mathcal H}}
\newcommand{\Z}{{\mathbb Z}}

\newcommand{\M}{{\mathcal M}}
\newcommand{\V}{{\mathcal V}}

\newcommand{\D}{{\mathbb D}}
\newcommand{\Dp}{{\overline{\Delta}_p}}
\newcommand{\IT}{{\mathbb T}}
\newcommand{\sphere}{{\widehat{\mathbb C}}}

\newcommand{\hull}{{\operatorname{hull}}}

\newcommand{\inv}{^{-1}}

\newcommand{\arccot}{{\,\operatorname{arccot}\,}}
\newcommand{\arth}{{\operatorname{arth}\,}}

\newcommand{\aand}{{\quad\text{and}\quad}}

\def\be{\begin{equation}}
	\def\ee{\end{equation}}

\newcommand{\bee}{\begin{enumerate}}
	\newcommand{\eee}{\end{enumerate}}

\newcommand{\blem}{\begin{lem}}
	\newcommand{\elem}{\end{lem}}
\newcommand{\bthm}{\begin{thm}}
	\newcommand{\ethm}{\end{thm}}
\newcommand{\bcor}{\begin{cor}}
	\newcommand{\ecor}{\end{cor}}
\newcommand{\beg}{\begin{examp}}
	\newcommand{\eeg}{\end{examp}}
\newcommand{\begs}{\begin{examples}}
	\newcommand{\eegs}{\end{examples}}
\newcommand{\bdefe}{\begin{defin}}
	\newcommand{\edefe}{\end{defin}}
\newcommand{\bprob}{\begin{prob}}
	\newcommand{\eprob}{\end{prob}}
\newcommand{\bei}{\begin{itemize}}
	\newcommand{\eei}{\end{itemize}}

\newcommand{\bcon}{\begin{conj}}
	\newcommand{\econ}{\end{conj}}
\newcommand{\bcons}{\begin{conjs}}
	\newcommand{\econs}{\end{conjs}}
\newcommand{\bprop}{\begin{propo}}
	\newcommand{\eprop}{\end{propo}}
\newcommand{\br}{\begin{rem}}
	\newcommand{\er}{\end{rem}}
\newcommand{\brs}{\begin{rems}}
	\newcommand{\ers}{\end{rems}}
\newcommand{\bo}{\begin{obser}}
	\newcommand{\eo}{\end{obser}}
\newcommand{\bos}{\begin{obsers}}
	\newcommand{\eos}{\end{obsers}}
\newcommand{\bpf}{\begin{pf}}
	\newcommand{\epf}{\end{pf}}
\newcommand{\ba}{\begin{array}}
	\newcommand{\ea}{\end{array}}
\newcommand{\beq}{\begin{eqnarray}}
	\newcommand{\beqq}{\begin{eqnarray*}}
		\newcommand{\eeq}{\end{eqnarray}}
	\newcommand{\eeqq}{\end{eqnarray*}}

\newcounter{minutes}
\divide\time by 60
\newcounter{hours}
\multiply\time by 60 \addtocounter{minutes}{-\time}
\begin{document}
	\title{Gehring-Hayman inequality for meromorphic univalent mappings}
	\begin{center}
		{\tiny \texttt{FILE:~\jobname .tex,
				printed: \number\year-\number\month-\number\day,
				\thehours.\ifnum\theminutes<10{0}\fi\theminutes}
		}
	\end{center}
	\author{Bappaditya Bhowmik${}^{~\mathbf{*}}$}
	\address{Bappaditya Bhowmik, Department of Mathematics,
		Indian Institute of Technology Kharagpur, Kharagpur - 721302,
		India.} \email{bappaditya@maths.iitkgp.ac.in}
	\author{Deblina Maity}
	\address{Deblina Maity, Department of Mathematics,
		Indian Institute of Technology Kharagpur, Kharagpur - 721302,
		India.}  \email{deblinamaity1997@gmail.com}
	\author{Toshiyuki Sugawa}
	\address{Toshiyuki Sugawa, Graduate School of Information Sciences, Tohoku University, Aoba-ku, Sendai - 9808579, Japan. } \email{sugawa@tohoku.ac.jp
	}
	
	%${}^{\mathbf{*}}$ Corresponding author}
%\date{ %\today Feb. 5, 08
	%June 23, 2016; File: MP-V1.tex}
%\date{%\today
	%September 01, 06; File: lau$_{-}$revise1.tex}
\begin{abstract}
	Let $f$ be a meromorphic univalent function on the open unit disk having a simple pole at $p\in (0,1)$ 
	that extends continuously to the left half $\IT^{-}$ of the unit circle.
	In this article, we prove that the ratio of the length of the image of the vertical diameter $\IA$ of the unit disk
	to the length of the image of $\IT^{-}$ under the mapping $f$ is bounded by a constant depending only on $p.$
	Next, we extend this result by considering any hyperbolic geodesic and any Jordan curve in $\D$ sharing the same endpoints. 
	These results extend the classical Gehring-Hayman inequality to meromorphic univalent functions
	and also prove a conjecture posed by Bhowmik and Maity [Bull. Sci. Math. \textbf{199} (2025), \# 103583].
\end{abstract}
%\thanks{The first author would like to thank SERB, India for its financial support through Core Research Grant
%	(Ref. No.- CRG/2022/001835)}
\thanks{The second author of this article would like to thank
	the Prime Minister Research Fellowship of the Government of India (Grant ID: 2403440) for its financial support.}

\maketitle
\pagestyle{myheadings} \markboth{B. Bhowmik, D. Maity, and T. Sugawa}%
{Gehring-Hayman inequality for meromorphic univalent mappings}

% \bigskip
%\noindent
\section{Introduction}
Let $\C, \D, \uhp$ and $\sphere$ stand for the complex plane, the open unit disk
$\{\zeta\in\IC: |\zeta|<1\},$ the upper half-plane $\{z\in\C: \Im z>0\}$
and the Riemann sphere $\C\cup\{\infty\},$ respectively, throughout this paper.
Also, let $\IT^{-}:=\{\zeta\in \C:~|\zeta|=1, ~\Re \zeta<0\}$ denote the left half 
of the unit circle $\IT=\partial\D$ and $\IA$ the diameter $(-i,i)$ of the unit circle with endpoints $i,-i$.
In other words,
$$
\IT^-=\{-e^{i\theta}: -\pi/2<\theta<\pi/2\}
\aand
\IA=\{iv: -1<v<1\}.
$$
We denote by $\ell(\gamma)$ the Euclidean length of a curve $\gamma$ in $\IC.$
Therefore, if $\gamma:(a,b)\to\C$ is a piecewise $C^1$ curve, then
$$
\ell(\gamma)=\int_a^b|\gamma'(t)|dt.
$$
In this paper, a conformal map $f$ of a domain $D$ means a holomorphic and
injective map of $D$ into $\sphere.$
A theorem of Gehring and Hayman \cite{GH} can be stated in the following form.
(Note that the original form deals with functions on the upper half-plane instead of the unit disk.)

\vspace{0.2cm}
\noindent\textbf{Theorem A} (Gehring-Hayman \cite[Theorem 1]{GH}).
\textit{
	Let $f$ be a conformal map of $\D$ into $\C$
	that remains continuous up to $\IT^-.$
	Then the inequality
	$$
	\ell(f(\IA))\le A_1 \ell(f(\IT^-))
	$$
	holds, where $A_1$ is an absolute constant with $\pi \leq A_1 < 74$.
}
\vspace{0.2cm}

Later, Jaenisch \cite{J} improved the range of $A_1$ to $4.56\leq A_1 \leq 17.45$. 
%We now recall the definition of the hyperbolic metric in $\D$ in order to present the generalized Gehring-Hayman theorem. The hyperbolic metric in the unit disk $\D$ is given by 
%$$
%\lambda_{\D}(z_1,z_2)=\underset{C}{\min}\int_{C}\frac{|dz|}{1-|z|^2},~~z_1,z_2\in\D,
%$$
%where the minimum is taken over all paths $C$ in $\D$ from $z_1$ to $z_2$; the curve that attains this minimum is called the hyperbolic geodesic from $z_1$ to $z_2$ (see \cite[p. 6]{Pom}). Now we state the generalized Gehring-Hayman theorem.

The second theorem of Gehring and Hayman \cite{GH} can also be stated as follows.

\vspace{0.2cm}
\noindent
\textbf{Theorem B} \cite[Theorem 2]{GH}.
\textit{Let $f$ be a conformal mapping of $\D$ into $\C$.
	Suppose that $J$ is a Jordan arc joining two points $\zeta_1$ and $\zeta_2$ in $\D.$
	Then
	$$
	\ell(f(S))\le A_1\ell(f(J)),
	$$
	%	$$	\int_{S}|f'(z)||dz|\leq A\int_{J}|f'(z)||dz|,	$$
	where $S$ is the hyperbolic geodesic joining $\zeta_1$ to $\zeta_2$ in $\D$
	and $A_1$ is the constant appearing in Theorem A.
}
\vspace{0.2cm}

We refer to the articles \cite{BM, PP} for various applications of Theorem B. 
% Also, see \cite{HR} for the Gehring-Hayman inequality concerning quasihyperbolic geodesics.
In the paper \cite{BD} of the first and the second authors, 
we tried to extend Theorems A and B to
meromorphic univalent mappings in $\D.$
To formulate our results, we define quantities $A_p$ for $0<p<1.$
We denote by $\M_p$ the class of
those meromorphic univalent functions in $\D$ with a simple pole at $z=p\in(0,1)$
that extend continuously to the left half $\IT^-$ of the unit circle.
We consider the problem to look for a constant $A_p$ such that
\be\label{eq:Ap}
\ell(f(\IA))\leq A_p\;\ell(f(\IT^-))
\ee
for all functions $f\in\M_p.$
We denote by $\M_p^0$ the subset of $\M_p$ consisting of $f$ with $\ell(f(\IT^-))<+\infty.$
Then the above condition remains unchanged if we replace $\M_p$ by $\M_p^0.$
In what follows, let $A_p$ be the best possible constant in \eqref{eq:Ap}.
That is,
$$
A_p=\sup_{f\in\M_p^0}\frac{\ell(f(\IA))}{\ell(f(\IT^-))}.
$$

%Now, instead of taking the horizontal diameter $(-1,1)$ and $\partial \D^{+}$ as in the Gehring-Hayman's theorem, 
%we consider the vertical diameter $ I:=\{z\in\D: {\rm{Re}}\,z=0; ~|{\rm{Im}}\,z|<1\}$ and the left half of the unit circle, i.e.  $C':=\{ z\in \partial \D:~ {\rm{Re}}\,z<0\}$. Let $\|\cdot\|$ denote the Euclidean length of a curve. 
In \cite[Theorem~1]{BD}, it was shown indeed that $A_p<+\infty$ for $\sqrt 2-1<p<1.$
More precisely, the result is stated as follows.

\vspace{0.2cm}
\noindent\textbf{Theorem C} (Bhowmik and Maity \cite[Theorem 1]{BD}).
\textit{Let $p\in (\sqrt{2}-1,1)$. Then
	\be\label{p1Th1Eq2}
	\frac{(1+p)^2 \pi}{4p}\leq A_p\leq \displaystyle{\min_{q\in(1,\infty)}R_p (q)},
	\ee
	where
	\be\nonumber
	R_p(q)=\frac{(1+p^2)\log q}{2p}\cot^2
	\left(\frac14\left(\tan^{-1}\left(\frac{q-1}{q+1}\right)-\tan^{-1}\left(\frac{(1-p^2)(q-1)}{2p(q+1)}\right)\right)\right).
	\ee }
\vspace{0.2cm}

Note that $(1-p^2)/2p\ge 1$ when $0<p\le\sqrt2-1.$
In \cite[p.~13]{BD}, it was conjectured that $A_p$ is finite for all $p\in(0,1).$
In view of the proof of Theorem C, we can even extend the class $\M_p$ to the class
$\tilde\M_p$ of analytic univalent functions on the domain $\D\setminus \Dp$
which extends to $\IT^-$ continuously,
where $\Dp=\{\zeta: |\zeta-\alpha\inv|\le \sqrt{\alpha^{-2}-1}\}$
and $\alpha=2p/(1+p^2).$ Note that $p\in\partial \Dp$ (see the next section) and thus $\M_p\subset\tilde\M_p$.
We consider the quantity
$$
\tilde A_p=\sup_{f\in\tilde\M_p^0}\frac{\ell(f(\IA))}{\ell(f(\IT^-))},
$$
where $\tilde\M_p^0$ denotes the subclass of $\tilde\M_p$ consisting of functions
$f$ such that $\ell(f(\IT^-))<+\infty.$
Since $\M_p^0\subset \tilde\M_p^0,$ we have $A_p\le \tilde A_p$ for $0<p<1.$
In this article, we improve the above bound of $A_p$ and prove the above conjecture.
We indeed show the following.

\bthm\label{thm:main1}%\label{p2Th1}
The following inequalities hold for each $p\in(0,1):$
\be
\frac{(1+p)^2 \pi}{4p}\leq A_p\le \tilde A_p
<\frac{1+p^2}{p}\left(1+\sqrt2+\frac{20}{3p}\right)^2\log2.
\ee
\ethm

Actually, we will give a better but much more complicated bound in Theorem~\ref{thm:main2} below.
We mention here that the proof of Theorem~\ref{thm:main2} will be accomplished by employing 
a similar technique adopted by Gehring and Hayman in \cite[Theorem~1]{GH}. 
%Moreover, instead of taking the vertical diameter $I$ and the left semicircle $C'$ in Theorem~\ref{p2Th1}, we consider any hyperbolic geodesic and any Jordan curve in $\D$ with the same endpoints and prove a generalized version of Theorem~B in Theorem~\ref{p2Th2}. 

The structure of the paper is as follows. Section 2 is devoted to several lemmas,
which will be used to prove our main results.
In Section 3, Theorem~\ref{thm:main2} will be presented together with its proof.
Then Theorem~\ref{thm:main1} will be derived from it.
Finally, as an application of the proof of Theorem~\ref{thm:main2},
we will extend Theorem B to the meromorphic case (see Theorem~\ref{thm:main3}).

\section{Preliminary results}
We denote by $\lambda_D(w)|dw|$ the hyperbolic metric of a domain $D$ in $\C$ with constant
Gaussian curvature $-4.$
The induced hyperbolic distance on $D$ will be denoted by $d_D(w_1,w_2).$
We recall that the hyperbolic distance in the unit disk $\D$ is given by
$$
d_\D(\zeta,\omega)=\arth \left|\frac{\zeta-\omega}{1-\zeta \bar\omega}\right|,~ ~\zeta,\omega\in\D,
$$
where $\arth x=({1/2})\log(({1+x})/({1-x}))$ for $0\le x<1.$
Let $f$ be a meromorphic univalent mapping of $\D$ with a simple pole at $p\in(0,1)$.
Choose $\alpha\in(0,1)$ so that $p$ is the hyperbolic midpoint of $0$ and $\alpha.$
In other words,
\be\label{Eq0.25}
p
=\frac{\alpha}{1+\sqrt{1-{\alpha}^2}}
=\frac{1-\sqrt{1-{\alpha}^2}}{\alpha}~\text{ i.e. }~
\alpha=\frac{2p}{1+p^2}.
\ee 
%Note that $\alpha\in(0,1)$ if and only if $p\in(0,1)$. 
Let $\gamma$ be the hyperbolic geodesic in $\D$ which is symmetric about the real axis and passes through the point 
$p$. 
Since hyperbolic geodesics of $\D$ are circular arcs in $\D$ perpendicular to $\partial \D$, we obtain 
\be\label{Eq0.5}
\gamma = \{\zeta\in\D: |\zeta-\alpha\inv|=\sqrt{\alpha^{-2}-1} \}.
\ee 
The two endpoints of $\gamma$ on $\partial \D$ are given by $\alpha \pm i\sqrt{1-{\alpha}^2}$. 
Let $\Omega$ be the connected component of $\D\setminus\gamma$ containing the origin.
In other words, $\Omega=\D\setminus\Dp,$ where
$\Dp=\{\zeta: |\zeta-\alpha\inv|\le \sqrt{\alpha^{-2}-1}\}$ as is given in Section 1.

The function
$$
\phi_\alpha(\zeta)=\frac{\zeta(1-\alpha \zeta)}{\zeta-\alpha}
$$
maps $\Omega$ conformally onto $\D$ and satisfies $\phi_\alpha(0)=0, \phi_\alpha'(0)=-\alpha\inv$
(\textit{cf.} \cite[p.~85]{MM}).
Then the M\"obius transformation
\beq \label{Eq5}
g(\zeta)=i\,\frac{1+i\zeta}{1-i\zeta}
\eeq
maps $\D$ conformally onto the upper half-plane $\mathbb{H}$ with $g(0)=i$ and
\be\label{Eq5.25}
g(p)
%g({\alpha}/{(1+\sqrt{1-{\alpha}^2}}))
=-\alpha+i\sqrt{1-\alpha^2}.
\ee
Note that $g$ is involutive; that is, $z=g(\zeta)$ if and only if $\zeta=g(z).$
In particular, $g(\uhp)=\D.$
Also, we see
\beq\nonumber
&g(\IA)&=\{iy: y>0\}=:\Ip,\\\nonumber
&g(\IT^-)&=\{x\in\IR: x>0\}=:{\IR}^{+}.
\eeq
Observe that the above curve $\gamma$ is mapped by $g$ onto the semi-circle 
$$%\be\label{}
\sigma = \{z\in\uhp: |z+\alpha\inv|=\sqrt{\alpha^{-2}-1}\}.
$$%\ee
We note that $\sigma$ is a hyperbolic geodesic in $\mathbb{H}$ with endpoints at 
$-p$ and $-1/p.$
%${-\alpha}/({1+\sqrt{1-{\alpha}^2}})=-p$ and $-({1+\sqrt{1-{\alpha}^2}})/{\alpha}=-1/p$. 
Let $\Omega_1$ be the hyperbolic half-plane $g(\Omega)$ in $\uhp.$
That is, $\Omega_1$ is the connected component of 
$\uhp\setminus\sigma$ containing $i.$

We now recall the concept of harmonic measure.
Let $D$ be a domain in $\mathbb{C}$, and let $E$ be a Borel subset of 
the boundary $\partial {D}$ of $D$. 
The harmonic measure of $E$ on $D$ is defined 
as the solution $u$ of the generalized Dirichlet problem in $D$ with boundary data $1$
on $E$ and $0$ on $\partial {D}\setminus E$ (see \cite[Ch.~3]{LV})
and its value $u(z)$ is denoted by $\omega(z,E,D)$ for $z\in D.$
Indeed, the set function $E\to \omega(z,E,D)$ can be regarded as a Borel probability
measure on $\partial D$ for each $z\in D.$
Let $D_1$ and $D_2$ be two domains in $\C.$
If $D_1\subset D_2$ and a Borel set $E$ is contained in 
$\partial D_1 \cap\partial D_2$, then $\omega(z,E,D_1)\leq\omega(z,E,D_2)$ 
for $z\in D_1$. 
This property is known as Carleman's principle of domain extension 
for harmonic measure (see \cite[p. 86]{Pom}).

For a proper subdomain $D$ of $\C,$ we will denote by $\delta_D(z)$ the Euclidean distance
from a point $z\in D$ to the boundary $\partial D.$
We denote by $B(a,r)$ the closed disk $|z-a|\le r$ in $\C$
centered at $a$ with radius $r>0.$

The proof of our main results will be based on the following estimate provided by
Gehring and Hayman \cite[Lemma 2]{GH}.

\vspace{0.2cm}
\noindent\textbf{Lemma~A.}\label{L2}
\textit{Let $D$ be a simply connected proper subdomain of $\C$ 
	and suppose that $d>0$ and $w_0\in \C\setminus D$ are given.
	For a Borel set $\beta$ contained in $\partial D\cap B(w_0,d)$ the
	inequality
	$$
	\delta_D(w)\leq d\,\cot^2({\pi}\omega(w,\beta, D)/{4})
	$$
	holds for every $w\in D.$
}
\vspace{0.2cm}

Next, we wish to establish an all-important lemma that will serve as 
a key ingredient in proving Theorem \ref{thm:main1}.

\blem \label{p2L3}
Let $\beta=[a,b]$ and $\beta'=[ia,ib]$ 
for positive numbers $a,b$ with $a<b.$
Then 
\be\label{eq:lem1}
\frac{\arccot M_p(b/a)}{\pi}
\leq\omega(z,\beta,\Omega_1)\leq\, \frac{1}{2}, \quad
z\in \beta',
\ee
where $\Omega_1$ is the hyperbolic half-plane defined above and
\be\label{eq:Mp}
M_p(q)=\frac{q+1}{q-1}+\frac{(1-p^2)^2(1+q^2)}{2p(q-1)(4p\sqrt{q}+(1+q)(1+p^2))}.
\ee
\elem

\begin{figure}[h]
	\centering
	\includegraphics[width=0.85\textwidth]{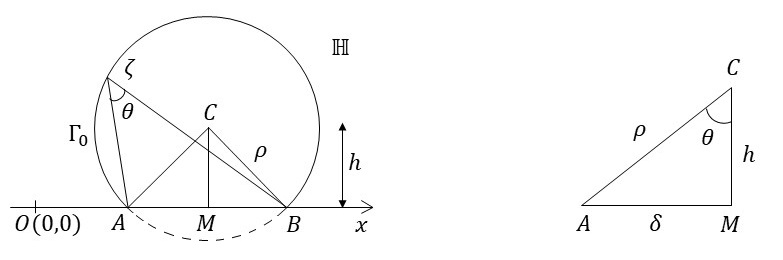} % Change 'image.jpg' to your file name
	\caption{The level curve $\Gamma_0$ and the triangle $\Delta CAM$}
	\label{fig_1}
\end{figure}

\bpf
Since $\Omega_1\subset\mathbb{H}$ and 
$\beta\subset \partial\Omega_1\cap\partial\uhp,$ 
Carleman's principle of domain extension for harmonic measure implies
\be\nonumber
\omega(z,\beta,\Omega_1)\leq\omega(z,\beta,\mathbb{H})
=\frac1\pi\arg\frac{b-z}{a-z}
\ee
for each $z\in\Omega_1.$
Note that $\omega(z,\beta,\mathbb{H})$ is ${1}/{\pi}$ times the angle 
subtended by the segment $\beta$ from the point $z.$
In particular, the set $\{z\in\uhp: \omega(z,\beta,\uhp)>1/2\}$ is equal to
the open upper half-disk with diameter $\beta=[a,b].$
Since $\beta'$ is away from this half-disk, we have
$\omega(z,\beta,\mathbb{H})\leq 1/2$ for $z\in \beta'$.
We now observe that the function
\be\nonumber
\psi(z)=\left(\frac{z+p}{z+1/p}\right)^2
\ee
maps $\Omega_1$ conformally onto $\mathbb{H}.$
Let $\tilde\beta=\psi(\beta)$ and $\tilde\beta'=\psi(\beta').$
Then the conformal invariance of harmonic measures implies
$$%\be \label{p2L3eq3}
\omega(z,\beta,\Omega_1)
%=\omega(\psi(z),\psi(\beta),\psi(\Omega_1))
=\omega(\psi(z),\tilde\beta,\mathbb{H}), \quad z\in \Omega_1.
$$%\ee
Therefore, it suffices to show the inequality $\pi\omega(\zeta,\tilde\beta,\uhp)
\ge\arccot M_p(b/a)$ on $\tilde\beta'$ for the lower estimate in \eqref{eq:lem1}.
To this end, we fix a point $\zeta$ in $\tilde\beta'$ and 
put $\theta=\pi \omega(\zeta,\tilde\beta,\uhp).$
Then the required inequality turns to the form $\cot\theta\le M_p(b/a).$
First, observe that the level curve
\be\nonumber
\Gamma_0=\{z\in\mathbb{H}: \omega(z,\tilde\beta,\mathbb{H})={\theta}/{\pi}\}
\ee
is an open subarc of the circle $\Gamma$ passing through the points
$A=\psi(a), B=\psi(b)$ and $\zeta.$
Let $C=M+ih$ be the center of $\Gamma$ as shown in Figure~\ref{fig_1}.
Note that $h>0$ because $\theta$ is an acute angle.
Now the inscribed angle theorem implies that
the line segment $\tilde\beta=[A,B]$ subtends the angle $2\theta$ at the center $C$.
Thus,
\be\label{p2L3eq4.5}
M=\frac{A+B}2 \aand
\tan \theta=\frac{B-A}{2h}.
\ee
We now compute the radius $\rho$ of $\Gamma.$
To this end, we put $\delta=(B-A)/2$
(see Figure~\ref{fig_1}) and compute
\be\nonumber
h^2=\rho^2-\delta^2
%=|\zeta-C|^2-\delta^2
=|\zeta-(M+ih)|^2-\delta^2
=|\zeta-M|^2-2h\Im(\zeta-M)+h^2-\delta^2.
\ee
By this relation, we find the value of $h$ as follows:
\be\nonumber
h=\frac{|\zeta-M|^2-\delta^2}{2\Im(\zeta-M)}
=\frac{|\zeta|^2-2\,M\Re\zeta+AB}{2\Im\zeta}.
\ee
Since \eqref{p2L3eq4.5} implies $\tan\theta=\delta/h,$ we now obtain
\be\label{eq:cot}
\cot \theta =\frac{|\zeta|^2-2\,M\Re\zeta+AB}{2\,\delta\Im\zeta}.
\ee
Here, a straightforward computation yields
$$
2M
=\psi(a)+\psi(b)
=\frac{p^2\Theta}{(1+ap)^2(1+bp)^2},
$$
where
$$
\Theta=(a^2+b^2)(1+p^4)+2(a+b)(1+ab)p(1+p^2)+2p^2(1+4ab+a^2b^2),
$$
$$
2\delta=\frac{(b-a)p^2(1-p^2)\big((a+b)(1+p^2)+2p(1+ab)\big)}%
{(1+ap)^2(1+bp)^2},
$$
and
$$
AB=\frac{(a+p)^2(b+p)^2 p^4}{(1+ap)^2(1+bp)^2}.
$$
Since $\zeta\in\tilde\beta'=\psi(\beta'),$ we can write $\zeta=\psi(iy)$
for some $y\in[a,b]$ so that
$$
\zeta=\frac{p^2(p^2+4p^2y^2+p^2y^4-y^2-p^4y^2)+
	2p^3(1-p^2)y(1+y^2)\,i}{(1+p^2y^2)^2}
$$
and
$$
|\zeta|=\frac{p^2(p^2+y^2)}{1+p^2y^2}.
$$
We substitute them into \eqref{eq:cot}, and after some simplification,
we obtain
$$
\frac{(b-a)\cot\theta}{a}
=\frac by+\frac ya+
\frac{(1-p^2)^2(a^2+y^2)(b^2+y^2)}{2ap\left(2p(1+ab)+(a+b)(1+p^2)\right)y(1+y^2)}.
$$
For convenience, we introduce two additional variables
$$
q=\frac ba, \quad
T=\frac ya\in[1,q]
$$

and compute

$$
(q-1)\cot\theta=\frac{(q+T^2)}{T}
+\frac{(1-p^2)^2(1+T^2)(q^2+T^2)}{2p Q},
$$

where

$$
Q=T({1}/{a}+aT^2)(2p({1}/{a}+aq)+(1+q)(1+p^2)).
%\frac{T(1+a^2T^2)(2p(1+a^2q)+a(1+q)(1+p^2))}{a^2}.
$$
Here we recall the elementary inequality 
$x+y\geq 2\sqrt{xy}$ for $x\geq0$ and $y\geq0$,
which will be used repeatedly below.
Then we get the lower estimate
$$
Q\ge 2T^2(4p\sqrt{q}+(1+q)(1+p^2))
$$
and hence
\begin{align*}
	(q-1)\cot\theta &\leq \frac{q+T^2}{T}
	+\frac{(1-p^2)^2(1+T^2)(q^2+T^2)}{4pT^2(4p\sqrt{q}+(1+q)(1+p^2))} \\
	%&=\frac{q}{T}+T+\frac{(1-p^2)^2({1}/{T}+T)({q^2}/{T}+T)}{4p(4p\sqrt{q}+(1+q)(1+p^2))} \\
	&= \frac{q}{T}+T+\frac{(1-p^2)^2(q^2/T^2+1+q^2+T^2)}{4p(4p\sqrt{q}+(1+q)(1+p^2))} \\
	&\le 1+q+\frac{(1-p^2)^2(1+q^2)}{2p(4p\sqrt{q}+(1+q)(1+p^2))}
\end{align*}
for $1\le T\le q.$
In the last part of the above estimations, we used the fact that the function
$c/x+x$ is convex for a constant $c>0$
and thus takes its maximum at one of the endpoints on the
closed interval $[x_0, x_1].$
Now the proof is complete.
\epf

We end this section by giving an elementary lemma, which will be utilized 
in the proof of Theorem \ref{thm:main1}.

\begin{lem}\label{lem:cot}
	Let $k$ be a constant with $0<k<1.$
	Then the function $G(x)=\cot(k \arccot x)$ is strictly convex and increasing 
	on $0<x<+\infty.$
	Moreover, $G(0^+)=\cot(k\pi/2), G'(0^+)=k/\sin^2(k\pi/2)$ and the following
	inequalities hold:
	$$
	\cot\frac{k\pi}2+\frac{kx}{\sin^2(k\pi/2)}<
	G(x)< \cot\frac{k\pi}2+\frac xk,\quad 0<x<+\infty.
	$$
\end{lem}

\bpf
First we compute $G(0^+)=\lim_{x\to0^+}G(x)=\cot(k\pi/2)$ because 
$\arccot(0^+)=\pi/2.$
Next we have
$$
G'(x)=\frac{k}{(1+x^2)\sin^2(k\arccot x)}>0,
$$
which implies that $G(x)$ is (strictly) increasing and gives the value of
$G'(0^+)$ in the assertion.
Since $\tan \theta$ is strictly convex on $0<\theta<\pi/2$ and satisfies $\tan 0=0,$
the inequality $\tan(k\theta)< k\tan\theta;$ equivalently, 
$\cot\theta< k\cot(k\theta),$ holds for $0<\theta<\pi/2.$
Letting $x=\cot \theta ,$ we obtain $x<k\cot(k\arccot x)$ for $0<x<+\infty.$
Hence, we observe that
$$
G''(x)=\frac{2k(k\cot(k\arccot x)-x)}{(1+x^2)^2\sin^2(k\arccot x)}>0
$$
for $0<x<+\infty.$
Thus we have checked strict convexity of $G.$
We next find the asymptotic value by using the L'H\^{o}pital's rule as follows:
\begin{align*}
	\lim_{x\to+\infty}\frac{G(x)}{x}
	&=\lim_{x\to+\infty}G'(x) 
	=\lim_{x\to+\infty}\frac{k}{(1+x^2)\sin^2(k\arccot x)} \\
	&=\lim_{\theta\to0^+}\frac{k}{(1+\cot^2\theta)\sin^2(k\theta)} 
	=\lim_{\theta\to0^+}\frac{k\sin^2\theta}{\sin^2(k\theta)}=\frac1k.
\end{align*}
The strict convexity of $G$ yields the inequalities
$$
G(0^+)+G'(0^+)(x-0)<
G(x)< G(0^+)+\frac{x(G(x_1)-G(0^+))}{x_1-0}
$$
for $0<x<x_1<+\infty.$
In particular, the left-hand inequality follows.
We now let $x_1\to+\infty$ to obtain the right-hand inequality.
\epf

\section{Proof of main results}
We now state and prove a technical result which leads to our main result.

\bthm\label{thm:main2}
Let $p\in (0,1).$ Then
\be\label{eq:Np}
\frac{(1+p)^2 \pi}{4p}\leq A_p\le \tilde A_p
\leq \displaystyle{\inf_{q\in(1,\infty)}N_p(q)},
\ee
where
$$
N_p(q)=\frac{(1+p^2)\log q}{2p} \cot^2\left(\frac14\cot^{-1}\left(\frac{q+1}{q-1}+
\frac{(1-p^2)^2(1+q^2)}{2p(q-1)(4p\sqrt{q}+(1+q)(1+p^2))}\right)\right).
$$
\ethm

\bpf
The proof of this theorem is similar to that of \cite[Theorem~1]{BD}. 
But, for the sake of completeness, we will repeat all the necessary details. 
Let $f\in\tilde\M_p$ for a given number $0<p<1.$
In other words, $f$ is a univalent analytic function on $\D\setminus\Dp$ 
%with $f(p)=\infty$
which is continuous up to the (open) left half $\IT^-$ of the unit circle.
Let $g:\D\to\uhp$ be the map defined in \eqref{Eq5} and let $\alpha$ be the number
defined in \eqref{Eq0.25}. 
Then $h=f\circ g^{-1}$ is analytic and univalent on $\Omega_1=g(\Omega)$ 
%with a pole at $\alpha'=-\alpha+i\sqrt{1-\alpha^2}$ (see \eqref{Eq5.25}) 
and remains continuous on ${\IR}^{+}$. 
%Assume $w_1=w|_{\Omega_1}$, where $\Omega_1$ is the same as in Lemma~A. 
%Since $\alpha'$ is contained in $\sigma,$
%the univalent function $h$ is bounded on $\Omega_1$ and continuous on ${\IR}^{+}$. 
In particular, the image $D=h(\Omega_1)=f(\Omega)$ is a simply connected domain in $\C.$
Choose $q\in(1,\infty)$ arbitrarily and fix it. 
Let $\beta_n=[q^n,q^{n+1}]$ and $\beta_n'=[iq^n,iq^{n+1}]$ for each integer $n.$
Let $\ell_n$ and $\ell'_n$ denote the lengths of the images of these segments 
under the mapping $h$. 
Then Lemma~\ref{p2L3} yields the inequalities
\be\label{eq:h}
\frac{\xi}{\pi}\le
\omega(h\inv(w),\beta_n,\Omega_1)=\omega(w,h(\beta_n), D)
\le 1/2
\ee
at each point $w\in h(\beta'_n)$, where $\xi=\cot^{-1}(M_p(q))$ 
and $M_p(q)$ is given in \eqref{eq:Mp}. 
%Note that $0<\xi<\pi/2$. 
Let $w_0$ be the point on $h(\beta_n)$ that divides its arclength in half. 
Then $h(\beta_n)$ lies in the disk $|w-w_0|\leq \ell_n/2$, and thus
Lemma A together with \eqref{eq:h} implies
\be\label{eq:delta}
\delta_{D}(w)\le 
\frac{\ell_n}{2}\cot^2(\pi\omega(w,h(\beta_n), D)/4)
\le \frac{\ell_n}{2}\cot^2(\xi/4)
\ee
for $w\in h(\beta'_n).$
We now recall that $\phi_\alpha$ maps $\Omega$ conformally onto the unit disk $\D$
and that $g=g\inv$ maps $\Omega$ onto $\Omega_1.$
Thus $\phi_\alpha\circ g:\Omega_1\to \D$ is conformal.
Since $h:\Omega_1\to D$ is conformal, too, we have the relations
\begin{align*}
	\lambda_D(h(z))|h'(z)|&=\lambda_{\Omega_1}(z)
	=\lambda_\D(\phi_\alpha(g(z)))|(\phi_\alpha\circ g)'(z)| \\
	&=\frac{2|\phi_\alpha'(w)|}{|z+i|^2(1-|\phi_\alpha(w)|^2)},
	\quad z\in \Omega_1,
\end{align*}
where we put $w=g(z).$
By the well-known estimate $1/\lambda_D\le 4\delta_D$ on $ D,$ we obtain
$$
\delta_D(h(z))\ge \frac1{4\lambda_D(h(z))}
=\frac{|h'(z)|(1-|\phi_\alpha(w)|^2)|z+i|^2}{8|\phi_\alpha'(w)|},
\quad
w=g(z)\in \Omega.
$$
We now take $z=iy\in\Ip$ and compute $w=g(iy)=i(1-y)/(1+y)=iv,$ where $v\in(-1,1).$
Hence, in conjunction with \eqref{eq:delta}, we get
$$
|h'(iy)|\le \frac{\ell_n}{2}\cot^2(\xi/4)
\frac{8|\phi_\alpha'(iv)|}{(1-|\phi_\alpha(iv)|^2)(y+1)^2}.
$$
Noting $dv/dy=-2/(y+1)^2,$ we estimate
\begin{align*}
	\ell'_{n}
	&=\int_{\beta'_{n}}|h'(z)||dz|
	=\int_{q^n}^{q^{n+1}}|h'(iy)|dy \\
	&\le\frac{\ell_n}{2}\cot^2(\xi/4)
	\int_{q^n}^{q^{n+1}}\frac{8|\phi_\alpha'(iv)|}{(1-|\phi_\alpha(iv)|^2)(y+1)^2}dy \\
	&=2\ell_n\cot^2(\xi/4)
	\int^{v_n}_{v_{n+1}}\frac{|\phi_\alpha'(iv)|}{1-|\phi_\alpha(iv)|^2}dv,
\end{align*}
where we put $v_n=(1-q^n)/(1+q^n)$ for $n\in\Z.$
Using the expression of $\phi_{\alpha},$ we have the inequality
$$
\frac{|\phi_\alpha'(iv)|}{1-|\phi_\alpha(iv)|^2}
=\frac{\sqrt{(1-v^2)^2+4\alpha^2 v^2}}{\alpha(1-v^4)}
\le \frac{\sqrt{(1-v^2)^2+4v^2}}{\alpha(1-v^4)}
=\frac{1}{\alpha(1-v^2)}
$$
for $v\in(-1,1)$ and $0<\alpha<1.$
Hence, using \eqref{Eq0.25}, we obtain
$$
\ell'_{n}\le 
2\ell_n\cot^2(\xi/4)
\int^{v_n}_{v_{n+1}}\frac{dv}{\alpha(1-v^2)}
=\ell_{n} \cot^2(\xi/4)\frac{\log q}{\alpha}=N_p(q)\ell_n
$$
for each $n\in\Z.$
Summing up for $n\in\Z,$ we finally get
$$
\ell(f(\IA))=\ell(h(\Ip))=\sum_{n=-\infty}^{\infty}\ell_n'
\le N_p(q)\sum_{n=-\infty}^{\infty}\ell_n
=N_p(q)\ell(h(\R^+))
=N_p(q)\ell(f(\IT^-)).
$$
The assertion of the theorem is now clear.
\epf

We now deduce Theorem \ref{thm:main1} from Theorem \ref{thm:main2}.

\bpf[Proof of Theorem \ref{thm:main1}]
By Theorem \ref{thm:main2} with the choice $q=4,$ we obtain
$A_p\le N_p(4)$ for $0<p<1,$ where
$$
N_p(4)=\frac{(1+p^2)\log 2}{p} \cot^2\left(\frac14\arccot
\left(\frac 53+\frac{17(1-p^2)^2}{6p(8p+5(1+p^2))}\right)\right).
$$
We now apply Lemma \ref{lem:cot} with $k=1/4$ to get
$$
N_p(4)<\frac{(1+p^2)\log 2}{p}\left(\cot\frac\pi 8+
4\left(\frac 53+\frac{17(1-p^2)^2}{6p(8p+5(1+p^2))}\right)\right)^2.
$$
Consider the function
$$
H(p)=6p\left(\frac 53+\frac{17(1-p^2)^2}{6p(8p+5(1+p^2))}\right)
=\frac{17p^4+50p^3+46p^2+50p+17}{5p^2+8p+5}.
$$
Since
$$
H(1)-H(p)=\frac{(1-p)(17p^3+67p^2+63p+33)}{5p^2+8p+5}>0
$$
for $0<p<1,$ we obtain $H(p)<H(1)=10.$
Noting the formula $\cot(\pi/8)=1+\sqrt 2,$ we now get
$$
N_p(4)<\frac{(1+p^2)\log 2}{p}\left(\cot\frac\pi 8+
\frac{2H(p)}{3p}\right)^2
<\frac{(1+p^2)\log 2}{p}\left(1+\sqrt 2+\frac{20}{3p}\right)^2.
$$
\epf

Table~\ref{table1} presents a comparison between ranges of $A_p$ by listing the upper bounds (UB) and lower bounds (LB) for selected values of $p \in (0,1)$ with the previously established results from \cite[Table 1]{BD}. As evident from the following table, the newly obtained upper bounds for $A_p$ show an improvement.

\begin{table}[ht]
	
	\begin{tabular}{|c |c
			|c |c |c |}		
		\hline 
		$p$ & LB of $A_p$  & \shortstack{UB of $A_p$\\ from Theorem~C} & \shortstack{UB of $A_p$ \\ from Theorem~1}  & \shortstack{UB of $A_p$\\ from Theorem~2} \\ [0.1ex] 
		\hline
		0.999 & 3.141 &  73.421 & 114.486  & 73.251\\
		\hline
		0.99 & 3.141 &  74.995 & 116.025 & 73.259\\
		\hline
		0.9 & 3.150 & 95.491 & 134.471 & 74.212 \\
		\hline
		0.8 & 3.180 & 135.733 & 164.134 & 77.634 \\
		\hline
		0.7 & 3.242 & 221.807 & 210.271 &  84.837 \\
		\hline
		0.6 & 3.351 & 471.016 & 287.415 &  98.455\\
		\hline
		0.5 & 3.534 & 1984.431 & 429.726 &  124.383\\
		\hline
		0.4 & 3.848 & ---& 731.847 & 178.045\\
		\hline
		0.3 & 4.424 &--- & 1528.574 & 310.577 \\
		\hline
		0.2 & 5.654 &--- & 4605.973 &775.275 \\
		\hline
		0.1 & 9.503 & --- & 33408.930 & 4608.760\\
		%\hline
		%0.01 & 80.118 & --- & 31033169 &  3535612\\
		\hline
	\end{tabular} 
	\vspace{0.2cm}\caption{Comparison between ranges of $A_p$ for various values of $p\in (0,1)$}\label{table1}
\end{table}

It is plausible that the quantity $A_p$ is decreasing in $0<p<1$ but it seems
difficult to show it rigorously.
Instead, we will show that $\tilde A_p$ is monotone.

\begin{lem}\label{lem:mono}
	The quantity $\tilde A_p$ is non-increasing in $0<p<1.$
	In particular, the quantity
	$$
	A_p^+=\sup_{p\le p'<1}A_{p'}
	$$
	satisfies the inequalities
	$$
	A_p\le A_p^+\le \tilde A_p
	$$
	for $0<p<1.$
\end{lem}

\bpf
It is enough to show the inclusion relation $\Delta_{p_1}\supset\Delta_{p_2}$
for $0<p_1<p_2<1.$
Here, we recall that $\Dp=\{\zeta: |\zeta-\alpha\inv|\le \sqrt{\alpha^{-2}-1}\}$
and $\alpha=2p/(1+p^2)$ for $0<p<1.$
To this end, we need to check the inequality
$\alpha_1\inv-\alpha_2\inv+\sqrt{\alpha_2^{-2}-1}\le \sqrt{\alpha_1^{-2}-1}$
for $\alpha_k=2p_k/(1+p_k^2),~ k=1,2.$
Indeed, it follows from the fact that the function $x-\sqrt{x^2-1}$ is decreasing
in $1<x<+\infty.$
\epf

\brs
\begin{enumerate}
	\item
	It is clear that $A_1\le \tilde A_p$ for each $0<p<1,$ where $A_1$ is the best possible
	constant appearing in Theorem A.
	In particular, we have from Theorem~\ref{thm:main2}
	$$
	A_1\le \tilde A_{1^-}\le \lim_{p\to 1^-}N_p(q)
	=\log q\cot^2\left(\frac14 \arccot\frac{q+1}{q-1}\right)
	$$
	for each $1<q<+\infty.$
	With the help of Mathematica, we choose $q=5.55465$ to obtain
	$\tilde A_{1^-}\le 73.2502105.$
	This agrees with the upper bound in Theorem~A.
	\item
	In Theorems~\ref{thm:main1} and \ref{thm:main2},
	we consider the configuration $\Dp$ and $\IA$ only.
	However, we can translate $\Dp$ by a disk automorphism of the form
	$T_a(z)=(z+ia)/(1-iaz)$ for $-1<a<1$ while leaving $\IA$ invariant.
	Hence, we will have the same result by replacing $\Dp$
	by $T_a(\Dp)$, where
	$$
	T_a(\Dp)=\left\{\zeta: \left|\zeta-\frac{1-a^2}{\alpha(1+a^2)}-\frac{2ai}{1+a^2}\right|\le \frac{\sqrt{\alpha^{-2}-1}(1-a^2)}{1+a^2}\right\}.
	$$
\end{enumerate}
\ers

Finally, we will generalize Theorem~B.
In order to enhance the difference,
we rephrase Theorem B in terms of the hyperbolic geometry of a simply connected
domain $D$ as follows.
Let $\Gamma$ be the hyperbolic geodesic that joins two points $w_1$ and $w_2$ in
a simply connected domain $D\subset\C.$
Then, for any Jordan arc $J$ joining $w_1$ and $w_2$ in $D,$ the inequality
\be\label{eq:geod}
\ell(\Gamma)\le A_1\ell(J)
\ee
holds, where $A_1$ is the constant in Theorem~A.
We will generalize this %the above form of Theorem A
to the case when $D$ is a simply connected hyperbolic subdomain of $\sphere$ containing $\infty.$
In this case, the bound should depend on a distance from $\infty$ to the geodesic $\Gamma.$
Our result can be stated in the form similar to that of Theorem B.
We denote by $\hull(C)$ for a closed curve $C$ in $\C$
the compact set $K=\C\setminus U_\infty,$ where $U_\infty$ is 
the (unique) unbounded connected component of $\C\setminus C,$
and denote by $\hull^\circ(C)$ the interior of the set $\hull(C).$
Here, we also note that the Euclidean length $\ell(J)$ of a Jordan arc $J$
is the same as the (suitably normalized) 1-dimensional Hausdorff measure $\mathcal{H}^1(J)$
(see \cite[p. 56]{Mat}).

\begin{thm}\label{thm:main3}
	Let $J$ be a (closed) Jordan arc in $\D$ with endpoints $\zeta_1$ and $\zeta_2$
	and let $\gamma$ be the (closed) hyperbolic line segment joining $\zeta_1$ and $\zeta_2$
	in $\D$.
	Let $\tilde\gamma$ be the set $C\setminus C_0,$ where $C$ is the Euclidean circle (or Euclidean line
	attached with $\infty$) containing $\gamma$ and $C_0$ is the connected component of $C\setminus J$
	containing a point outside $\D.$
	Let $f$ be a meromorphic univalent function on $\D$ with a pole
	at $s\in\D.$
	Suppose that $s\not\in\hull(\tilde\gamma\cup J)$,
	then
	$$
	\ell(f(\gamma))\le
	\begin{cases}
		A_1\ell(f(J)) & \text{if}~ s\not\in \hull^\circ(J\cup \hat J) \\
		A_\tau^+\ell(f(J)) & \text{if}~ s\in \hull^\circ(J\cup \hat J),
	\end{cases}
	$$
	where $\hat J$ denotes the reflection of $J$ with respect to the circle $C$ and
	$$
	\tau=\tanh d_\D(s,\tilde\gamma).
	$$
\end{thm}

Here, $A_1$ is the same as defined in Theorem~A and $A_p^+$ is defined in Lemma~\ref{lem:mono} for $0<p<1.$
If $A_p$ is non-increasing, we would have $A_p^+=A_p.$

\bpf
By our assumption, $\tau>0.$ Translating by a disk automorphism if necessary,
we can assume that $\zeta_1$ and $\zeta_2$ lie on $\IA=(-i,i).$
Then $\hat J$ is the reflection of $J$ in the imaginary axis.
Let $\V$ be the set of all connected components $V$ of $J\setminus\tilde\gamma.$
Then each $V\in\V$ is an open subarc of $J$ with endpoints in $\tilde\gamma.$
We denote by $I_V$ the open line segment sharing two endpoints with $V.$
Let $I=\cup_{V\in \V}I_V.$
Then $I$ is an open subset of $\tilde\gamma.$
We set $\gamma_0=\tilde\gamma\cap J$ and show that
$$
\gamma\subset\gamma_0\cup I \quad \mbox{ and }\quad 
J=\gamma_0\cup \left( \bigcup_{V\in\V} V\right).
$$
Since $\cup_{V\in\V}V=J\setminus\tilde\gamma,$ the latter equality is clear. We will verify the inclusion relation $\gamma\subset\gamma_0\cup I.$ 
Choose an arbitrary point $\zeta_0\in \gamma\setminus I$ and we will show the claim
$\zeta_0\in J.$
%We suppose, to the contrary, that $\zeta_0\not\in J.$
We denote by $\tilde\zeta_m~ (m=1,2)$ the endpoints of $\tilde\gamma$ so that
$\tilde\zeta_1, \zeta_1, \zeta_2, \tilde\zeta_2$ lie in $\IA$ in this order.
Note that $\tilde\zeta_m$ may be equal to $\zeta_m.$
Let $\tilde\gamma_m=[\tilde\zeta_m,\zeta_0]$ for $m=1,2.$
For convenience, we parametrize $J$ by a homeomorphism $j:[0,1]\to J$ such that
$j(0)=\zeta_1$ and $j(1)=\zeta_2.$
%Take $t_0\in[0,1]$ so that $j(t_0)=\zeta_0.$
Let $t_1=\sup\{t\in[0,1]: j(t)\in\tilde\gamma_1\},$
that is, $t_1$ is the ``final exit time" of $j$ from $\tilde\gamma_1.$
In particular, $j(t)\ne \zeta_0$ for $t>t_1.$
It suffices to show that $\zeta_0=j(t_1)~ (\in J).$
To the contrary, we suppose that $j(t_1)\ne \zeta_0.$
We next let $t_2=\inf\{t\in(t_1,1]: j(t)\in\tilde\gamma_2\}$.
By continuity of $j,$ we see that $\zeta_m':=j(t_m)\in \tilde\gamma_m$
and $j(t_m)\ne\zeta_0$ for $m=1,2.$
In particular, $t_1<t_2.$
Since $j((t_1,t_2))$ does not meet $\tilde\gamma$ by construction,
$V=j((t_1,t_2))\in\V$ and $I_V=(\zeta_1',\zeta_2')$ contains the point
$\zeta_0,$ which contradicts the assumption $\zeta_0\in\gamma\setminus I.$
Hence, we conclude that $\zeta_0=j(t_1)\in J.$
Thus, the claim has been verified and $\zeta_0\in \gamma\cap J\subset\gamma_0$ follows.

We now fix a Jordan arc $V\in\V.$
We denote by $D_V$ the Jordan domain bounded by $J_V=\overline{V}\cup I_V$ and denote by
$\widetilde D_V$ the union of $D_V,~ I_V$ and $\hat D_V$
which stands for the reflection of $D_V$ in $\IA.$
Let $g_V:\D_-\to D_V$ be a conformal homeomorphism whose continuous extension
to $\partial\D_-$ maps $\IA$ onto $I_V.$
Here, we recall that $\D_-$ is the left-half of the unit disk $\D.$
Schwarz's reflection principle now allows us to extend $g_V$ to a conformal
homeomorphism of $\D$ onto $\widetilde D_V$, which will be denoted 
by the same symbol.
We next consider the following two cases.

\noindent
{\bf Case I}: $s\not\in \widetilde D_V.$
In this case, $f\circ g_V:\D\to\C$ is analytic and univalent.
Thus, by Theorem A, we obtain
$$
\ell(f(I_V))=\ell((f\circ g_V)(\IA))\le A_1\ell((f\circ g_V)(\IT^-))=A_1\ell(f(V)).
$$
\noindent
{\bf Case II}: $s\in \widetilde D_V.$
In this case, since $s\not\in D_V\cup I_V$ by assumption,
we may further normalize the conformal map $g_V:\D\to\widetilde D_V$ 
so that $p=g_V\inv(s)\in(0,1).$
Since $\widetilde D_V\subset\D,$ the domain monotonicity of hyperbolic metric
yields
$$
\arth p=d_\D(p,0)=d_{\widetilde D_V}(s,g_V(0))\ge
d_{\widetilde D_V}(s,I_V)
\ge d_\D(s,I_V)\ge d_\D(s,\tilde\gamma)=\arth \tau,
$$
and therefore, $p\ge \tau.$
By the definition of $A_p$, we have
$$
\ell(f(I_V))=\ell((f\circ g_V)(\IA))\le A_p\ell((f\circ g_V)(\IT^-))=A_p\ell(f(V))
\le A_\tau^+\ell(f(V)).
$$

Finally, we show the inequality in the assertion.
When $s\in \hull^\circ(J\cup \hat J),$ we have
$\ell(f(I_V))\le A_\tau^+\ell(f(V))$ for every $V\in\V$
because $A_1\le A_\tau^+.$
Since $\gamma\subset\gamma_0\cup I$ and $J$ is the disjoint union
of $\gamma_0$ and $V$ over $\V,$ we obtain
\begin{align*}
	\ell(f(\gamma))&=\Hm^1(f(\gamma))\le\Hm^1(f(\gamma_0))+\Hm^1(f(I)) \\
	&\le\Hm^1(f(\gamma_0))+\sum_{V\in\V}\Hm^1(f(I_V)) \\
	&\le A_\tau^+\Hm^1(f(\gamma_0))+A_\tau^+\sum_{V\in\V}\Hm^1(f(V)) \\
	&=A_\tau^+ \Hm^1(f(J))=A_\tau^+\ell(f(J)).
\end{align*}
When $s\not\in \hull^\circ(J\cup \hat J),$ only Case I occurs.
Therefore, we can show the corresponding inequality in a similar way.
\epf

\noindent{\bf Statements and Declarations}: Nil.\\
\noindent{\bf Competing interests}: We have nothing to declare.\\
\noindent{\bf Conflict of interest}: Nil.\\

\end{document}